\documentstyle[12pt]{article}

\newcommand{\const}{\mathop{\rm const}\limits}

\newcommand{\mes}{\mathop{\rm mes}\limits}

\newcommand{\vraisup}{\mathop{\rm vraisup}\limits}

\newcommand{\card}{\mathop{\rm card}\limits}

\newcommand{\supp}{\mathop{\rm supp}\limits}

\newcommand{\Ent}{\mathop{\rm Ent}\limits}

\begin{document}

\begin{center}

{\bf FOURIER SERIES AND TRANSFORMS IN GRAND LEBESGUE SPACES} \\
\vspace{3mm}
{\bf As an particular case - exponential Orlicz spaces.}\\

\vspace{3mm}

 $ {\bf E.Ostrovsky^a, \ \ L.Sirota^b } $ \\

\vspace{4mm}

$ ^a $ Corresponding Author. Department of Mathematics and computer science, Bar-Ilan University, 84105, Ramat Gan, Israel.\\
\end{center}
E - mail: \ eugostrovsky@list.ru\\
\begin{center}
$ ^b $  Department of Mathematics and computer science. Bar-Ilan University,
84105, Ramat Gan, Israel.\\
\end{center}
E - mail: \ sirota@zahav.net.il\\

\vspace{3mm}

{\bf Abstract}.  In this article we investigate the Fourier series and transforms for
the functions defined on the $ [-\pi, \pi]^ d $ or on the $ R^d $ and  belonging to  the
(Bilateral) Grand Lebesgue Spaces.\par
 As a particular case we obtain some results about Fourier's transform in
 the so-called exponential Orlicz spaces. \par
 We construct also several examples to show the exactness of offered estimations.\par

\vspace{3mm}


\vspace{3mm}

\vspace{3mm}

\section{Introduction. Notations. Problem Statement.} \par
\vspace{3mm}
 For the and real valued measurable function $ f = f(x), \  $ defined on the
 $ X = \{ x \} = T^d =[-\pi,\pi]^d $ or equally $ X = [0, \ 2 \pi]^d,
\ d = 1,2, \ldots $ or $ X = R^d $ we denote correspondingly the Fourier coefficients and
transform
$$
c(n)=\int_{T^d} \exp( i(n,x) ) \ f(x) \ dx, \ \ F[f](t) = \int_{R^d} \exp(i(t,x)) \ f(x) \ dx,
$$
 where as usually
$$
F[f](t) \stackrel{def}{=}  \lim_{M \to \infty}\int_{|x| \le M} \exp(i(t,x)) \ f(x) \ dx,
$$

$$
x = (x_1,x_2,\ldots,x_d), \  n = (n_1,n_2, \ldots, n_d),  \ t = (t_1,t_2, \ldots, t_d), \
dx = \prod_{j=1}^d dx_j
$$
if $ x \in R^d, $ and $ dx = (2 \pi)^{-d} \ \prod_{j=1}^d dx_j $ in the case $ X = T^d; $
$$
n_i = 0, \pm 1, \pm 2, \ldots, \ (t,x) = \sum_{j=1}^d t_j x_j, \ |n| = \max_{j} |n_j|, \
|t| = \max_j |t_j|,
$$

$$
s_M[f](x) = (2 \pi)^{-d}  \sum_{n: |n| \le M} c(n) \exp(-i (n,x)), \ x \in T^d;
\eqno(1.1)
$$

$$
 S_M[f](x) = (2 \pi)^{-d} \int_{t: |t| \le M} \exp(-i(t,x)) \ F[f](t) \ dt, \ x \in R^d.
 \eqno(1.2)
$$

\vspace{2mm}
{\bf Our aim in this paper is investigating of
the boundedness of (linear)  operators $ S_M[\cdot], \ s_M[\cdot] $ in  some Grand
Lebesgue Spaces norms $ ||\cdot||G\psi $ (see definition further):

$$
\sup_{M \ge 2} \sup_{f: ||f||G\psi = 1} ||s_M[f]||G\psi_1 < \infty, \
\sup_{M \ge 2} \sup_{f: ||f||G\psi = 1} ||S_M[f]||G\psi_1 < \infty,\eqno(1.3)
$$
 and the convergence and divergence in this norms}
 $$
 s_M[f](\cdot) \to f, \ S_M[f](\cdot) \to f  \eqno(1.4)
 $$
 {\bf as } $ M \to \infty.$ \par

\vspace{2mm}

We will prove that the so-called exponential Orlicz spaces over $ X $ are the
particular cases of the Grand Lebesgue Spaces. Therefore,
we can  consider also that the function $ f(\cdot) $ belongs to some Orlicz space
$ L(N) = L(N;X) $   with so - called {\it exponential
$ N \ - $ function } $ N = N(u), $
and will investigate the properties of Fourier transform of $ f, $ for example,
the boundedness of operators $ S_M[\cdot], \ s_M[\cdot] $ and
the convergence and divergence (1)in {\it some} Orlicz norms $ L(N;X). $ \par

\vspace{3mm}
  Note than the case if the function $ N(\cdot) $ satisfies the $ \Delta_2 $ condition is
known; see, for example, \cite{Ryan1}, \cite{Rao2}. Our results are also some
generalization of \cite{Gord1}, \cite{Ostrovsky2}, \cite{Pinsky1}, \cite{Wolff1} etc.\par

\vspace{3mm}

 The papier is organized as follows. In the next section we recall used facts 
 about Grand Lebesgue Spaces and obtain some new properties of this spaces, especially, 
 investigate the properties of these spaces in the case when the measure is discrete.
In the third section we obtain the GLS boundedness of Hilbert's transform.\par
 The fourth section is devoted to the weight Fourier operators boundedness in GLS spaces. The $ 5^{th} $ section contain the main result of the offered papier: the boundedness 
 of Fourier transforms in GLS spaces in general case, for instance, in the exponential 
 Orlicz spaces. In the next section we formulate and prove some auxiliary facts.\par
  The $ 7^{th} $ section contain the proofs of main results. \par
 In the last section we prove the GLS boundedness of the so-called maximal Fourier 
 operators.\par
 In many offered estimations we show their exactness by means of construction of 
 suitable (counter) examples.\par

\section{Grand Lebesgue Spaces.}

\vspace{3mm}

 Now we will describe using Grand Lebesgue Spaces (GLS) and a particular case the
 so-called Exponential Orlicz  Spaces (EOS). \\

\vspace{2mm}

 {\bf 1. Description of used Classical Lebesgue Spaces.} \par

 \vspace{2mm}

  Let $ (X,A,\mu) $ be some measurable space with sigma-finite
non - trivial measure $ \mu. $ For the measurable real valued
function $ f(x), \ x \in X, f: X \to R $ the symbol
$ |f|_p = |f|_p(X,\mu) $ will denote the usually $ L_p $ norm:
$$
|f|_p = ||f||L_p(X, \mu) = \left[ \int_X |f(x)|^p \ \mu(dx) \right]^{1/p}, \ p \ge 1.
\eqno(2.1)
$$
In the case $ X = R^d $ we introduce a new measure $ \nu(\cdot) $ (non - finite, in
general case): for all Borel set $ A \subset R^d $
$$
\nu(A) = \int_A \prod_{i=1}^d x_i^{-2} \ dx = \int_A \prod_{i=1}^d x_i^{-2} \ \cdot
\prod_{i=1}^d dx_i, \eqno(2.2)
$$
and will denote  $ |f|_p(\nu) = $
$$
\left| \left| \left(\prod_{i=1}^d x_i \right) \cdot f \right| \right|L_p(X,\nu) =
\left[ \int_X \left|\prod_{j=1}^d x_j \right|^p \ \cdot |f(x)|^p \ \nu(dx)  \right]^{1/p} =
$$
$$
\left[\int_X \left| \prod_{j=1}^d x_j \right|^{p-2} \cdot |f(x)|^p  \ dx    \right]^{1/p}. \eqno(2.3)
$$
 For arbitrary multiply sequence (complex, in general case)
$ c(n) = c(n_1,n_2,\ldots, n_d), n_i = 0, \pm 1, \pm 2, \ldots,
n \in Z^d $ we denote as usually
$$
|c|_p  =  \left[ \sum_{n} |c(n) |^p \right]^{1/p}, \ p \ge 1;\eqno(2.4)
$$
and introduce the discrete analog of $ |f|_p(\nu) $ norm:
$$
|c|_{p,\nu} = |c|^{(d)}_{p,\nu} =
\left[ \sum_n \left|c(n) \right|^p \ \cdot  \left(
\left| \prod_{j=1}^d n_j \right|^{p-2} +1 \right) \right]^{1/p}, \ p \ge 2.\eqno(2.5)
$$

\vspace{2mm}

{\bf 2. Grand Lebesgue Spaces.}\par

\vspace{2mm}
We recall in this section  for reader conventions
some definitions and facts from the theory of GLS spaces.\par
\vspace{2mm}

     Recently, see \cite{Fiorenza1},
     \cite{Fiorenza2}, \cite{Fiorenza3}, \cite{Iwaniec1}, \cite{Iwaniec2},
     \cite{Kozachenko1},
     \cite{Ostrovsky1}, \cite{Ostrovsky2}, \cite{Ostrovsky3}, \cite{Ostrovsky4},
     \cite{Ostrovsky5}, \cite{Ostrovsky6}, \cite{Ostrovsky9}  etc.
     appears the so-called Grand Lebesgue Spaces $ GLS = G(\psi) =G\psi =
    G(\psi; A,B), \ A,B = \const, A \ge 1, A < B \le \infty, $ spaces consisting
    on all the measurable functions $ f: T \to R $ with finite norms

     $$
     ||f||G(\psi) \stackrel{def}{=} \sup_{p \in (A,B)} \left[ |f|_p /\psi(p) \right].
     \eqno(2.6)
     $$

      Here $ \psi(\cdot) $ is some continuous positive on the {\it open} interval
    $ (A,B) $ function such that

     $$
     \inf_{p \in (A,B)} \psi(p) > 0, \ \psi(p) = \infty, \ p \notin (A,B).
     $$
We will denote
$$
 \supp (\psi) \stackrel{def}{=} (A,B) = \{p: \psi(p) < \infty, \} \eqno(2.7)
$$

The set of all $ \psi $  functions with support $ \supp (\psi)= (A,B) $ will be
denoted by $ \Psi(A,B). $ \par
  This spaces are rearrangement invariant, see \cite{Bennet1}, and
  are used, for example, in the theory of probability \cite{Talagrand1}, \cite{Kozachenko1}, \cite{Ostrovsky1}; theory of Partial Differential Equations \cite{Fiorenza2}, \cite{Iwaniec2}; functional analysis \cite{Ostrovsky4}, \cite{Ostrovsky5}; theory of Fourier series \cite{Ostrovsky7}, theory of martingales \cite{Ostrovsky2} etc.\par
 Notice that in the case when $ \psi(\cdot) \in \Psi(A,B),  $ a function
 $ p \to p \cdot \log \psi(p) $ is convex, and  $ B = \infty, $ then the space
$ G\psi $ coincides with some {\it exponential} Orlicz space. \par
 Conversely, if $ B < \infty, $ then the space $ G\psi(A,B) $ does  not coincides with
 the classical rearrangement invariant spaces: Orlicz, Lorentz, Marzinkievitch etc.

 We will use the following two important examples (more exact, the {\it two families
of examples} of the $ \psi $ functions and correspondingly the GLS spaces.\par
{\bf 1.} We denote
$$
\psi(A,B; \alpha,\beta;p) \stackrel{def}{=} (p-A)^{-\alpha} \ (B - p)^{-\beta},\eqno(2.8)
$$
where $ \alpha,\beta = \const \ge 0, 1 \le A < B < \infty; p \in (A,B) $ so that
$$
\supp \psi(A,B; \alpha,\beta;\cdot) = (A,B).
$$

{\bf 2.} Second example:
$$
\psi(1,\infty; 0, -\beta;p) \stackrel{def}{=} p^{\beta}, \eqno(2.9)
$$
but here $ \beta = \const > 0, \ p \in (1,\infty) $ so that
$$
\supp \psi(1,\infty; 0,-\beta;\cdot) = (1,\infty).
$$
 The space $ G\psi(1,\infty; 0, -\beta;\cdot)  $ coincides up to norm equivalence with
 the Orlicz space over the set  $ D $ with usually Lebesgue measure and with the
 correspondent $ N(\cdot) $ function

 $$
 N(u) = \exp\left(u^{1/\beta} \right), \ u \ge 1; N(u) = C |u|, |u| \le 1.
 $$
 Recall that the domain $ D $ has finite measure; therefore the behavior of the function
$ N(\cdot) $ is'nt essential.\par
{\bf Remark 1.} If we define the {\it degenerate } $ \psi_r(p), r = \const \ge 1 $ function as follows:
$$
\psi_r(p) = \infty, \ p \ne r; \psi_r(r) = 1
$$
and agree $ C/\infty = 0, C = \const > 0, $ then the $ G\psi_r(\cdot) $ space coincides
with the classical Lebesgue space $ L_r. $ \par
{\bf Remark 2.} Let $ \xi: D \to R $ be some (measurable) function from the set
$ L(p_1, p_2), \ 1 \le p_1 < p_2 \le \infty. $ We can introduce the so-called
{\it natural} choice $ \psi_{\xi}(p)$  as as follows:

$$
\psi_{\xi}(p) \stackrel{def}{=} |\xi|_p; \ p \in (p_1,p_2).
$$

\vspace{2mm}

{\bf 3. Discrete Grand Lebesgue Spaces.}\par

\vspace{2mm}

{\sc A. General part.}\par

\vspace{2mm}

 Let $ c = \vec{c} = \{ c(1), c(2), c(3), \ldots, c(n), \ldots  \} $ be arbitrary
 numerical sequence,
 $ \beta = \vec{\beta} = \{ \beta(1), \beta(2), \beta(3), \ldots, \beta(n), \ldots  \} $
 be arbitrary non-negative non-trivial:
 $$
 \sum_{n=1}^{\infty} \beta(n) \in (0,\infty]
 $$
 numerical sequence, $ p \in (A,B), \ 1 \le A < B
 \le \infty, \ \psi: (A,B) \to R_+, \ \psi \in \Psi(A,B). $ We define as before the
 so-called {\it weight discrete} GLS space $ G_d\psi_{\beta}(A,B) = G_d\psi_{\beta} $
 as a set of numerical sequences with finite norm

$$
|| c||G_d\psi_{\beta} = \sup_{p \in (A,B)}  \left[\frac{|c|_{p,\beta}}{\psi(p)} \right],
\eqno(2.10)
$$
where

$$
|c|_{p,\beta} \stackrel{def}{=}
\left[ \sum_{n=1}^{\infty} |c(n)|^p \ \beta(n)  \right]^{1/p}.\eqno(2.11)
$$

 Evidently, the $ G_d\psi_{\beta} $ spaces are particular cases of general GLS spaces,
 relative the weight measure

 $$
 \mu_{\beta}(A) = \sum_{k \in A} \beta(k).
 $$

 But this spaces are {\it resonant spaces } in the terminology of the book \cite{Bennet1} only in the case when $ \beta(n) = \const > 0. $ We can suppose in this case
without loss of generality that $ \beta(n) = 1 $ and will write for simplicity

$$
|| c||G_d\psi = \sup_{p \in (A,B)}  \left[\frac{|c|_p}{\psi(p)} \right],
$$
where as ordinary

$$
|c|_p \stackrel{def}{=}
\left[ \sum_{n=1}^{\infty} |c(n)|^p \right]^{1/p}.
$$

\vspace{2mm}

{\sc B. Natural function.}\par

\vspace{2mm}

 Let $ c = \vec{c} $ be the numerical sequence such that for some number
 $$
\exists  p_0 \in [1,\infty) \ |c|_{p_0} < \infty. \eqno(2.12)
 $$
 We investigate in this pilcrow the natural function
 $  \psi_c(p) $ for the sequence $ \vec{c}:$

 $$
 \psi_c(p) = |c|_p = \left[ \sum_{n=1}^{\infty} |c(n)|^p \right]^{1/p},
 $$
in addition to the assertions of the pilcrow 2. \par
 Note first of all that if $ q > p \ge 1, $ then $ |c|_q \le |c|_p. $ Therefore, if
for some  $ p_0 \in [1,\infty) \ |c|_{p_0} < \infty, $ then for all the values
$ p, \ p > p_0 \ \Rightarrow \psi_c(p) < \infty. $ \par
 Further,

$$
\lim_{n \to \infty} c(n) = 0,
$$
and
$$
\lim_{p \to \infty} |c|_p = \sup_n |c(n)| = \max_n |c(n)| \stackrel{def}{=} |c|_{\infty}.
$$

Thus, we proved the following assertion.\par

{\bf Lemma 1.} Every non-trivial natural discrete function $ \psi = \psi(p) $ has the
following properties: \\
0. The domain of definition of the  function $ \psi(\cdot) $ is some semi-axis
$ (p_0, \infty) $ or $ [p_0, \infty), $ where $ p_0 \ge 1. $ \\
1. The function $ \psi(\cdot) $ is monotonically non-increasing.\\

\vspace{2mm}
The proposition of the Lemma 1 is false in the case of weighted discrete GLS spaces. Let
us consider the correspondent example.\par
\vspace{2mm}

{\bf Example 1.} Let us consider the following weight sequence $ \beta^{(s)}: $

$$
\beta^{(s)}(n) = n^{-1 - s},
$$
and the following numerical sequence $ y = \{y(n)\}, $ where

$$
y(n) = n^{\theta}.
$$
Here $ s, \theta = \const, p_0 \stackrel{def}{=} s/\theta > 1.  $ \par

 Note that the norm $ ||y||_{p,\beta^{(s)} }, \ p \ge 1 $  is finite only when
$ p < p_0: $

$$
||y||^p_{p,\beta^{(s)}} = \sum_{n=1}^{\infty} n^{-1-s+p \theta}< \infty \ \Leftrightarrow p < s/\theta = p_0.
$$

{\bf  Example 2.} Let us consider the following sequence:
$ a = a^{(L)} = \{a(n)\}, a(n) = n^{-1/L}, \ L = \const \ge 1. $ We have for the values
$ p > L, $ denoting

$$
\psi_{(L)}(p) = \psi_{a^{(L)}}(p):
$$

$$
\psi^p_{(L)}(p) = \sum_{n=1}^{\infty} n^{-p/L}.
$$

 The last expression coincides with the well-known Rieman's zeta-function at the point
 $ p/L. $ Therefore,

 $$
 \psi_{(L)}(p) \asymp \left[\frac{p}{p-L} \right]^{1/L}, \ p \in (L,\infty).
 \eqno(2.13)
 $$

\vspace{2mm}

{\sc C. Tail behavior.}\par

\vspace{2mm}

 Let $ c = \vec{c}  \in G_d \psi(a,b), \ 1 \le a < b \le \infty. $ We introduce as ordinary the tail function $ T_{\beta}(c,u), u \in (0,\infty) $ for the sequence
 $ \{c\} $ relative an arbitrary discrete measure $ \mu_{\beta}: $
$$
 T_{\beta}(c,u) \stackrel{def}{=} \mu_{\beta}(n: |c(n)| \ge u).\eqno(2.14)
$$
 We will write for simplicity in the case $ \beta(n) = \beta_0(n) = 1, n = 1,2,\ldots $

$$
 T_{\beta_0}(c,u) \stackrel{def}{=} T(c,u) = \card \{ n: |c(n)| > u \}.
$$

 If the sequence $ \{ |c(n)| \} $ is bounded, for example, if for some $ p \ge 1 \
 |c|_p < \infty, $

 $$
 \forall u >  \sup_n |c(n)| \ \Rightarrow T(c,u) = 0.
 $$
 Therefore, we must investigate {\it in this case} the asymptotical behavior of the tail function $ T(c,\epsilon) $ only as $ \epsilon \to 0+. $ \par
 It follows from Tchebychev's inequality that
$$
T_{\beta}(c,u) \le \inf_{p \in (a,b)}
\left[ ||f||^p_{p,\beta} \psi^p(p)/u^p  \right],  \ u > 0. \eqno(2.15)
$$
 Conversely,
$$
|c|^p_{p,\beta} = p \int_0^{\infty} u^{p-1} T_{\beta}(c,u) du, \ p \ge 1;
$$
 therefore
$$
||c||G_d(\psi) =  \sup_{p: \psi(p)< \infty}
\left[ p \ \left[\int_0^{\infty} u^{p-1}
T_{\beta}(c,u) \ du \right]^{1/p} \ /\psi(p) \right]. \eqno(2.16)
$$
 In order to show the exactness of this inequalities, we consider some examples. \par
{\sc Example 1.} Let
$$
 a(k) =k, \ \beta(k)=k^{-b-1}, \ b = \const > 1.
$$
 We find by the direct calculations:

 $$
 |a|^p_{p,\beta} \sim (b-p)^{-1}, \ p \to b-0;
 $$

$$
T_{\beta}(a,u) = \sum_{k \ge u} k^{-b-1} \sim u^{-b}/b, \ u \to \infty;
$$
but it follows from the upper estimation for the tail function that as $ u \to \infty $

$$
T_{\beta}(a,u) \le C \inf_{p \in (1,b)} \left[ u^{-p} /(b-p) \right] \sim
C_1 u^{-b} \ \log u.
$$

 More generally, if
 $$
 \beta(k) = \beta^{(\Delta)}(k) = k^{-b-1} \ \log^{\Delta}(k),
 $$
$ b = \const > 1, \Delta = \const \ge 0, \ a(k) = k, $ and  $ p \to b - 0, $  then

$$
|a|^p_{p,\beta} = \sum_{k=1}^{\infty} k^{-b-1+p} \log^{\Delta}k \sim
$$
$$
\int_1^{\infty} x^{p-b-1} \ \log^{\Delta}(x) \ dx =
\frac{\Gamma(\Delta+1)}{(b-p)^{\gamma}}, \ \gamma \stackrel{def}{=}\Delta+1;\eqno(2.17)
$$

$$
T_{\beta^{(\Delta)}}(a,u) = \sum_{k \ge u} k^{-b-1} \log^{\Delta}(k) \sim
$$

$$
b^{-1} u^{-b} \ \log^{\Delta} u = b^{-1} u^{-b} \ \log^{\gamma-1} u, \ u \to \infty;
\eqno(2.18)
$$
but the upper estimation for the tail function gives only the inequality

$$
T_{\beta^{(\Delta)}}(a,u) \le C_2(b,\gamma) \ u^{-b} \ \log^{\gamma} u, \ u \ge e.
\eqno(2.19)
$$
 Let us show now that the inequality (2.18) is asymptotically exact as $ u \to \infty, $
 by virtue of the consideration of a following example. \par
Let us denote
 $$
 X(k) = \Ent \left[e^{e^k} \right],
 $$
where $ \Ent[z] $ denotes the integer part of the variable $ z;$
$$
p(k):= \exp \left( \gamma b k - b \exp(k) \right),
$$
where
$$
\gamma = \const > 0, \ b = \const > 1.
$$

 We define the weight sequence $ \beta(k) $ as follows:
 $$
 \beta(X(k)) = p(k)
 $$
and
$$
\beta(l) = 0, \ l \ne X(k) \ \forall k = 1,2,\ldots.
$$

 We introduce also the sequence $ y(n) = n, n = 1,2,\ldots. $  It is easy to compute
 analogously to \cite{Ostrovsky2}:

 $$
 |y|_{p, \beta} = \left\{\sum_k k^p \beta(k) \right\}^{1/p} \asymp (b-p)^{-\gamma},
  \ p \in (1,b),
 $$
but we observe that for the {\it subsequence} $ X(k) $
$$
T_{\beta}(y,X(k))= T_{\beta}(y,X(k)) \ge C_1(\gamma,b) \ [\log X(k)]^{\gamma}
\cdot X(k)^{-b}. \eqno(2.20)
$$

Note that the "continuous case" was investigated in
 \cite{Ostrovsky2}, \cite{Ostrovsky9}. \par
\vspace{2mm}
{\sc Example 2.} We know that for the sequence $ a(n) = n^{-1/L}, L = \const > 1 $

$$
|a|^p_p = \psi^p_{a^{(L)}}(p) \sim \frac{p}{p-L}, \ p > L;
$$
therefore we obtain from the upper tail estimation

$$
T(a,\epsilon) \le C \ \epsilon^{-L} \ |\log \epsilon|, \ \epsilon \in (0,1/e);
$$
but really

$$
T(a,\epsilon) \asymp  \epsilon^{-L}, \ \epsilon \to 0+.
$$
 In the more general case when the sequence $ a(n) $ has a view

 $$
 a(n) = n^{-1/L} \ \log^q n,  L > 1, q > 0,
 $$
we have:

$$
|a|^p_p \sim  \left[\frac{p}{p-L} \right]^{pq+1},\ p > L;
$$

$$
T(a,\epsilon) \sim \epsilon^{-L} \ |\log \epsilon|^{qL}, \epsilon \to 0+;
$$
Note that it follows from the upper estimation only the inequality

$$
T(a,\epsilon) \le C  \epsilon^{-L} \ |\log \epsilon|^{1+qL}, \epsilon \in (0, 1/e).
$$

 But we can show  that the our upper bound for the tail function is non-improvable.
 Namely, in the article \cite{Ostrovsky9} was constructed for all the values
 $ L = \const > 1, \ q \ge 0 $ the  example of {\it discrete}
  function $ z = z(k), \ k=1,2,\ldots $  and the correspondent weight
  $ \beta = \beta(k),$  for which

$$
|z|^p_{p,\beta} = \psi^p_{\beta}(z,p) \sim C_3(L,q) \
\left[\frac{p}{p-L} \right]^{pq+1}, \ p > L \eqno(2.21)
$$
and simultaneously for some positive {\it subsequence} $ \epsilon(m) $ monotonically
tending to zero

$$
T_{\beta}(z,\epsilon(m)) \ge C(L,q) \  \epsilon(m)^{-L} \ |\log \epsilon(m)|^{1+qL}, \epsilon(m) \in (0, 1/e). \eqno(2.22)
$$

 \vspace{2mm}

{\sc D. Leindler's inequality for discrete GLS spaces.}\par

\vspace{2mm}

 Let $ \beta = \{\beta(n) \}, n = 1,2,\ldots $ be again a discrete weight. We introduce
 a two linear operators:

 $$
 T[x](n) = \sum_{k=n}^{\infty} \frac{x(k) \beta(k)}{\Sigma(k)},
 $$
 where

 $$
 \Sigma(k) = \sum_{j=1}^k \beta(j);
 $$

 $$
 U[x](n) = \sum_{k=1}^{n} \frac{x(k) \beta(k)}{\sigma(k)},
 $$
where

$$
\sigma(k) = \sum_{i=k}^{\infty} \beta(i),
$$
in the case when

$$
 \sum_{i=1}^{\infty} \beta(i) < \infty.
$$

 Suppose that for some function $ \psi \in \Psi $
 $$
 x \in G_d(\psi);
 $$
for instance $ \psi $ may be the natural function for the sequence $ x(\cdot):
 \ \psi(p) = \psi_x(p). $ \par
 Denote
 $$
 \psi_1(p) = p \psi(p).
 $$

{\bf  Theorem L.} (Leindler's inequality for discrete GLS spaces.) \par
{\bf A.}

$$
||T[x]||G_d(\psi_1) \le 1 \cdot ||x||G_d(\psi), \eqno(2.23)
$$
where the constant "1" is the best possible.\\

{\bf B.}

$$
||U[x]||G_d(\psi_1) \le 1 \cdot ||x||G_d(\psi), \eqno(2.24)
$$
where the constant "1" is the best possible.\par
 Note that Leindler's inequalities for discrete GLS spaces are used for obtaining the
 $ L_p$ weight estimations for trigonometric  series, see  \cite{Yu1}. \par

{\bf Proof of the upper estimate. }\par
 0.  We will use the Leindler  inequalities \cite{Leindler1}, (which are some generalizations of the classical Hardy-Littlewood inequalities):

$$
\sum_{n=1}^{\infty} \beta(n) \left( \sum_{k=1}^n \alpha(k) \right)^p \le p^p \cdot \sum_{n=1}^{\infty} \beta^{1-p}(n) \alpha^p(n) \left(\sum_{k=n}^{\infty} \beta(k)  \right)^p; \eqno(2.25)
$$

$$
\sum_{n=1}^{\infty} \beta(n) \left( \sum_{k=n}^{\infty} \alpha(k) \right)^p \le p^p \cdot \sum_{n=1}^{\infty} \beta^{1-p}(n) \alpha^p(n) \left(\sum_{k=1}^n \beta(k)  \right)^p;
\eqno(2.26)
$$
In this inequalities $ \alpha(n) $ is arbitrary non-negative sequence, $ 1 \le p < \infty.$ \par
1. Let us prove the assertion {\bf A} of our theorem; the second may be proved analogously. Note that we can assume that all the variables $ x(n) $ are non-negative.\par
 We substitute in (2.25)

$$
\alpha(n) = \frac{\beta(n)}{\Sigma(n)} \cdot x(n),
$$
where
$$
x(\cdot) \in G_d(\psi),
$$
as long as in other case is nothing to prove.\par
We can and will suppose without loss of generality that $ ||x||G_d(\psi)\le 1, $
or equally

$$
\forall p \ \Rightarrow |x|_{p,\beta} \le \psi(p).
$$

The right-hand side $ R_a^{(p)} $ of inequality(2.25) has a view:
$$
R_a^{(p)} = p^p \ ||x||^p_{p,\beta},
$$
but the left-hand side $ L_a^{(p)} $ of this inequality may be rewritten as follows:

$$
L_a^{(p)} =  ||T[x]||^p_{p,\beta}.
$$
 We conclude using the first Leindler's inequality

 $$
 ||T[x]||_{p,\beta} \le p \cdot ||x||_{p,\beta} = p \ \psi(p) = \psi_1(p),
 $$
and after  using the direct definition of the norm in GLS spaces,

$$
||T[x]||G_d\psi_1 \le ||x||G_d\psi.
$$

\vspace{2mm}

{\bf Proof of the exactness.} \par
 We describe here the method of the lower estimations which will be used
 often further.\par

Let us denote

$$
V(x;\beta,\psi) = \frac{||T[x]||G_d(\psi_1)}{||x||G_d(\psi)},\eqno(2.27)
$$

$$
\overline{V} = \sup_{x \ne 0} \sup_{\psi \in \Psi(1,\infty)} \sup_{\beta > 0}
V(x;\beta,\psi),\eqno(2.28)
$$
and analogously

$$
V_0(x;\beta,\psi) = \frac{||U[x]||G_d(\psi_1)}{||x||G_d(\psi)},\eqno(2.29)
$$

$$
\overline{V}_0 = \sup_{x \ne 0} \sup_{\psi \in \Psi(1,\infty)} \sup_{\beta > 0}
V_0(x;\beta,\psi).\eqno(2.30)
$$
 From theorem L follows that

 $$
 V \le 1;  \ V_0 \le 1.
 $$

It remains to prove the inverse inequalities.\par
 Note first of all that the expression for the value $ V(\cdot) $ may be rewritten
 as follows:

$$
V(x;\beta,\psi) = \frac{\sup_p \left[|T[x]|_{p,\beta}/(p \psi(p)) \right]}
{\sup_p \left[ |x|_{p,\beta}/\psi(p) \right]}, \eqno(2.31)
$$
and if we choose

$$
\psi(p) = |x|_{p,\beta},
$$
i.e. $ \psi(\cdot) $ is the natural function for the sequence $ x $ relative the
weight $ \beta(\cdot): \ \psi = \psi_x, $  we obtain the following lower estimation for the value $ \overline{V}.$\par

{\bf Proposition 1.}

$$
\overline{V} \ge \sup_p W(p), \eqno(2.32)
$$
where the functional $ W =W(p) = W(p; T) $  has (here) a view:
$$
W(p) = \sup_{x \ne 0} \sup_{\beta > 0}
\left[  \frac{|T[x]|_{p,\beta}}{p |x|_{p,\beta}} \right].\eqno(2.33)
$$

 As a consequence: let $ x_0  $ be arbitrary element of the space $ l_{p,\beta} $ and
 $ \beta_0  $ be any sequence satisfying  our conditions, then $ W(p) \ge W_0(p), $
 where
 $$
W_0(p) = \left[  \frac{|T[x_0]|_{p,\beta_0}}{p |x_0|_{p,\beta_0}} \right].
\eqno(2.34)
$$
  Furthermore, if $ x_0^{\Delta}  $ be arbitrary {\it set:} $ \Delta = \const $ of
 elements of the space $ l_{p,\beta} $ and
$ \beta_0^{\Delta}  $ be any set of the
sequences satisfying  our conditions, then $ W(p) \ge W_1(p), $  where

 $$
W_1(p) = \sup_{\Delta}
\left[ \frac{|T[x_0^{\Delta}]|_{p,\beta_0^{\Delta}}}{p|x_0^{\Delta}|_{p,\beta_0^{\Delta}}} \right] \eqno(2.35)
$$
and consequently

$$
W_1(p) \ge W_2(p) \stackrel{def}{=} \overline{\lim}_{\Delta \to \infty}
\left[ \frac{|T[x_0^{\Delta}]|_{p,\beta_0^{\Delta}}}
{p |x_0^{\Delta}|_{p,\beta_0^{\Delta}}} \right] \eqno(2.36)
$$

Note in addition if there is some value $ p_0, \ p_0 \in \supp \psi  $  or the point
$ p_0 = \infty, $ in the case when $ \supp \psi = (A,\infty), $  which will be called
{\it critical point,} then

$$
\overline{V} \ge \overline{\lim}_{p \to p_0 \pm 0} W(p) \ge
\overline{\lim}_{p \to p_0 \pm 0} W_1(p) \ge
 \overline{\lim}_{p \to p_0 \pm 0} W_2(p). \eqno(2.37)
$$

\vspace{3mm}

{\bf We return to the proof of assertion of the considered theorem L. }\par
Taking here as an examples the values
$$
\beta(n) = n^s, \ \alpha(n)= n^{-1-\theta},
$$
where $ s, \theta = \const > 0, $ and

$$
p_0 \stackrel{def}{=} (s+1)/\theta > 1,
$$
we obtain after simple calculations:

$$
\sigma(n) = \sum_{k=n}^{\infty} k^{-1-s} \sim n^{-s}/s, \ n \to \infty;
$$

$$
x(n) = \alpha(n) \sigma(n)/\beta(n) \sim n^{1+\theta}/s,
$$
and we have as $ p \to p_0-0: $
$$
|x|^p_{p,\beta} \sim s^{-p} \sum_{n=1}^{\infty} n^{-1-s + p(1+\theta)} \sim
\frac{s^{-p}}{s-p(1+\theta)}.
$$
Further,
$$
T[x](n)= \sum_{k=1}^n k^{\theta} \sim n^{\theta+1}/(\theta+1), \ n \to \infty;
$$
$ p \to p_0 - 0 \ \Rightarrow $

$$
|T[x]|^p_{p,\beta} \sim \sum_{n=1}^{\infty} (\theta + 1)^{-p}n^{-1-s+p(\theta+1)} \sim
\frac{(\theta+1)^{-p}}{s-p(1+\theta)}.
$$
Thus,

$$
\overline{V} \ge \overline{\lim}_{p \to p_0 - 0}
\frac{s}{\theta+1} \cdot \frac{1}{p}=\lim_{p \to p_0-0} \frac{p_0}{p} = 1.
$$
\vspace{2mm}
 Analogously may be proved the estimate $ \overline{V}_0 \ge 1. $  It is sufficient to choose

 $$
 \beta(n) = n^t, \ \alpha(n) = n^{-1-\tau}, \ p_0 = (t+1)/\tau >1,
 $$
and
$$
\overline{V}_0 \ge  \overline{\lim}_{p \to p_0+0}\frac{t+1}{\tau} \cdot \frac{1}{p}= 1.
$$

\vspace{2mm}

{\bf 4. Exponential Orlicz Spaces.}\par

\vspace{2mm}
 We will prove in this subsection that the so-called Exponential Orlicz Spaces
 (EOS) are particular cases of Grand Lebesgue Spaces. \par
  In the case of finite measurable spaces, for example, for the probabilistic spaces
  this assertion was proved in \cite{Kozachenko1}; see also \cite{Ostrovsky1},
  chapter 1, section 5.\par
   Let $ N = N(u) $ be some
$ N \ - $ Orlicz's function, i.e. downward convex, even, continuous
differentiable for all sufficiently greatest  values $ u, \ u \ge u_0, $
strongly  increasing in the right - side axis,
and such that $ N(u) = 0 \ \Leftrightarrow  u = 0; \ u \to \infty \ \ \Rightarrow
dN(u)/du \ \to \ \infty. $ We say that  $ N(\cdot) $ is an Exponential Orlicz Function,
briefly: $ N(\cdot) \in \ EOF, $ if $ N(u) $ has a view: for some continuous
differentiable strongly increasing downward convex in the domain $ [2,\infty] $
function $ W = W(u) $ such that $ u \to \infty \ \Rightarrow W^/(u) \to \infty $

$$
N(u) = N(W,u) = \exp(W(\log |u|)), \ |u| \ge e^2. \eqno(2.38)
$$
 For the values $ u \in [-e^2,e^2] $ we define $ N(W,u) $ arbitrary  but so that the function
$ N(W,u) $ is even continuous convex strictly increasing in the right side axis
and such that $ N(u) = 0 \
\Leftrightarrow u = 0. $ The correspondent Orlicz space on $ T^d, \ R^d $ with usually
Lebesgue measure  with $ N \ - $ function $ N(W,u) $
we will denote $ L(N) = EOS(W); \ EOS = \cup_{W} \{EOS(W)\} $ (Exponential Orlicz's Space).\par
 For example, let $ m = \const > 0, \ r = \const \in R^1, $
$$
N_{m,r}(u) = \exp \left[|u|^m \ \left(\log^{-mr}(C_1(r) + |u| ) \right) \right] -1,
\eqno(2.39)
$$
$ C_1(r) = e, \ r \le 0; \ C_1(r) = \exp(r), \ r > 0. $ Then $ N_{m,r}(\cdot) \in EOS. $
In the case $ r = 0 $ we will write $ N_m = N_{m,0}. $ \par
 Recall here that the Orlicz's norm on the arbitrary measurable space $ (X,A,\mu) $
$ ||f||L(N) = L(N,X,\mu) $ may be calculated by the formula (see, for example,
\cite{Krasnoselsky1},p. 73; \cite{Rao1}, p. 66)
$$
||f||L(N) = \inf_{v > 0}  \left \{ v^{-1} \left(1 + \int_X N( v|f(x)|) \ \mu(dx) \right) \right \}. \eqno(2.40)
$$
 Recall also that the notation $ N_1(\cdot) << N_2(\cdot) $ for two Orlicz functions
$ N_1, N_2 $ denotes:
$$
\forall \lambda > 0 \ \Rightarrow \lim_{u \to \infty} N_1(\lambda u)/N_2(u) = 0.
\eqno(2.41)
$$
 We will denote for arbitrary Orlicz $ L(N) $ (and other r.i.) spaces by $ L^0(N) $ the
closure of all bounded functions with bounded support. \par
 Let $ \alpha $ be arbitrary number, $ \alpha = \const \ge 1, $ and $ N(\cdot)
\in EOS(W) $ for some $ W = W(\cdot). $ We denote for such a function
$ N = N(W,u) $
by $ N^{(\alpha)}(u) $ a new $ N \ - $ Orlicz's function such that
$$
N^{(\alpha)}(u) = C_1 \ |u|^{\alpha}, \ \ |u| \in [0, C_2];
$$
$$
N^{(\alpha)}(u) = C_3 + C_4 |u|, \ \ |u| \in (C_2, C_5];
$$
$$
N^{(\alpha) }(u) = N(u), \ \ |u| > C_5, \ \ 0 < C_2 < C_5 < \infty,\eqno(2.42)
$$
$$
C_{1,2,3,4,5} = C_{1,2,3,4,5}(\alpha,N(\cdot)).
$$
 In the case $ \alpha = m(j+1), \ m > 0, \ j = 0,1,2, \ldots $ the function $ N^{(\alpha)}_m(u) $
is equivalent to the following Trudinger's function:
$$
N_m^{(\alpha)}(u) \sim N_{[m]}^{(\alpha)}(u) = \exp \left(|u|^m \right) -\sum_{l=0}^{j}
u^{ml}/l!. \eqno(2.43)
$$
 This method is described in \cite{Rao1},
 p. 42 - 47. Those Orlicz spaces are applied to the theory
of non - linear partial differential equations.\par
 We can define formally the spaces $ L(N^{(\alpha)}_m) $ at $ m = + \infty $ as a projective
limit at $ m \to \infty $ the spaces $ L(N^{(\alpha)}_m), $ but it is evident that
$$
L \left(N^{(\alpha)}_{\infty} \right) \sim  L_{\alpha} + L_{\infty},
$$
where the space $ L_{\infty} $ consists on all the a.e. bounded functions with norm
$$
|f|_{\infty} = \vraisup_{x \in X} |f(x)|.
$$
 Of course, in the case $ X = T^d $
$$
L_{\alpha} + L_{\infty} \sim L_{\infty}.
$$

\vspace{3mm}

{\sc Hereafter we will denote by $  C_k = C_k(\cdot), k = 1,2,\ldots $ some positive finite essentially and by $ C, C_0 $ non-essentially "constructive" constants.\par
  By the symbols $ K_j = K_j(d) $ we will denote the "classical" {\it absolute} constants; more exactly, positive finite functions depending only on the dimension} $ d. $ \par

\vspace{3mm}

 It is very simple to prove the existence of constants $ C_{1,2,3,4,5} = C_{1,2,3,4,5}(\alpha,N(\cdot)) $ such that $ N^{(\alpha)} $ is some new exponential
 $ N \ $ Orlicz's function. \par

  Now we will introduce some {\it new } Grand Lebesgue Spaces. Let $ \psi = \psi(p),
\ p \ge \alpha, \alpha = \const \ge 1 $
be some continuous positive: $ \psi(\alpha) > 0 $  finite strictly increasing function such that the function $ p \to p \log \psi(p) $ is downward convex and
$$
\lim_{p \to \infty} \psi(p) = \infty.
$$
The set of all those functions we will denote $ \Psi; \ \Psi = \{\psi\}. $ A particular case:

$$
\psi(p) = \psi(W;p) = \exp(W^*(p)/p),
$$
where
$$
W^*(p) = \sup_{z \ge \alpha} (pz - W(z)) \eqno(2.44)
$$
is so - called Young - Fenchel, or Legendre transform of $ W(\cdot). $ It follows from
theorem of Fenchel - Moraux that in this case
$$
W(p) = \left[ p \ \log \psi(W;p) \right]^*, \ \ p \ge p_0 = \const \ge 2,
\eqno(2.45)
$$
and consequently  for all $ \psi(\cdot) \in \Psi $ we introduce the correspondent
$ N \ - $ function by equality:
$$
N([\psi]) = N([\psi],u) =
\exp \left \{\left[p \log \psi(p) \right]^*(\log u) \right \}, \ u \ge e^2.\eqno(2.46)
$$
Since $ \ \forall \ \psi(\cdot) \in \Psi, \ d = 0,1,\ldots \ \Rightarrow p^d \cdot \psi(p) \in \Psi, $ we can denote
$$
\psi_d(p) = p^d \cdot \psi(p), \ \  N_d([\psi]) = N_d([\psi],u) = N([\psi_d], u).
$$
For instance, if $ N(u)= \exp(|u|^m), \ u \ge 2, $ where $ \ m = \const > 0, $ then
$$
N_d([\psi], u) \sim \exp \left(|u|^{m/(dm+1)} \right), \ u \ge 2.
$$
{\bf Definition.} We introduce for arbitrary such a function $ \psi(\cdot) \in \Psi $
the so - called $ G(\alpha;\psi) \ $ and $ G(\alpha; \psi,\nu) $
 norms  and correspondent Banach spaces $ G(\alpha; \psi), \ G(\alpha, \psi, \nu) $ as a
set of all measurable (complex) functions with finite norms:

$$
||f||G(\alpha; \psi) = \sup_{p \ge \alpha} (|f|_p/\psi(p)),\eqno(2.47)
$$
and  analogously
$$
||f||G(\alpha; \psi,\nu) = \sup_{p \ge \alpha} (|f|_p(\nu)/\psi(p)).\eqno(2.48)
$$
 For instance $ \psi(p) $ may be $ \psi(p) = \psi_m(p) = p^{1/m}, \ m = \const > 0; $ in
this case we will write $ G(\alpha, \psi_m) = G(\alpha,m) $ and
$$
||f||G(\alpha,m) = \sup_{p \ge \alpha}  \left(|f|_p \ p^{-1/m} \right).\eqno(2.49)
$$
 Also formally  we define
$$
||f||G(\alpha, m) = |f|_{\alpha} + |f|_{\infty}.
$$
{\it Remark 1.} It follows from I'ensen inequality that in the case $ X = T^d $ all the spaces
$ G(\alpha_1; \psi), \ G(\alpha_2,\psi), 1 \le \alpha_1 < \alpha_2 < \infty $ are isomorphic:
$$
||f||G(\alpha; \psi) \le ||f||G(1; \psi) \le \max(1, \psi(\alpha)) \
||f||G(\alpha; \psi). \eqno(2.50)
$$
 It is false in the case $ X = R^d.$ \par
{\it Remark 2. } $ G(\alpha; \psi) $ is a rearrangement invariant (r.i.) space. $ G(\alpha,m) $ has a  fundamental function $ \phi(\delta; G(\alpha,m)), \ \delta > 0, $ where for any rearrangement invariant space $ G $
$$
\phi(\delta;G) \stackrel{def}{=} ||I_A(\cdot)||G(\cdot), \ \  \mes(A) = \delta \in (0,\infty),
$$
$ \mes(A) $ denotes usually Lebesgue measure of Borel set $ A. $  We have:
$$
 \phi(\delta; G\psi) = \sup_p [\delta/\psi(p)], \delta \in (0,\infty).
$$
 The detail investigating of $ G\psi \ $  spaces, for instance, their fundamental
 functions see in \cite{Liflyand1}, \cite{Ostrovsky2}.\par

 Let us consider also another space
$ G(a,b,\alpha, \beta), \ 1 \le a < b < \infty; \ \alpha,\beta \ge 0. $ Here
$ X = R^d $ and we denote $ h = \min((a+b)/2; 2a). $  We introduce the function
$ \zeta: (a,b) \to R^1_+: $
$$
 \zeta(p)= \zeta(a,b,\alpha,\beta; p) = (p-a)^{\alpha}, \ p \in (a, h);
$$

$$
 \zeta(p) = (b-p)^{\beta}, \ p \in [h, b).
$$
 By definition, the space $ G(a,b,\alpha,\beta) $ consists
on all the measurable complex functions with finite norm:
$$
||f||G(a,b,\alpha,\beta) = \sup_{p \in (a,b)} \left[|f|_p \ \cdot
\zeta(a,b,\alpha, \beta; p) \right].
$$
 The space $ G(a,b,\alpha,\beta) $ is also a rearrangement invariant space. \par

 For example, let us consider the function $ f(x) = f(a,b; x), \ x \in R^1 \to R: $
$ f(x) = 0, \ x \le 0; $
$$
f(x) = x^{-1/b}, \ x \in (0,1); \ \ f(x) = x^{-1/a}, \ x \in [1,\infty);
$$
then $ f(a,b,\cdot) \in G(a,b,1,1) $ and
$$
\forall \ \Delta \in (0, 1/2] \ \Rightarrow f \notin G(a,b,1-\Delta,1) \cup G(a,b, 1, 1-\Delta).
$$

 Analogously may be defined the "discrete" $ g(a,b,\alpha,\beta) $ spaces. Namely, let
$ c = c(n) = c(n_1,n_2,\ldots,n_d) $ be arbitrary multiply (complex) sequence. We say that $ c \in g(a,b,\alpha,\beta) $  if

$$
||c||g(a,b, \alpha, \beta) \stackrel{def}{=} \sup_{p \in (a,b)}
\left[ |c|_p \ (p-a)^{\alpha} \ (b-p)^{\beta} \right].
$$
 It is evident that the non - trivial case of those spaces is only if $ \beta = 0; $
in this case we will write $ g(a,b, \alpha, 0) = g(a, \alpha) $ and
$$
||c||g(\alpha) = \sup_{p > a} |c|_p \ (p-a)^{\alpha}.
$$
 We denote also for $ \psi(\cdot) \in \Psi: \ ||c||g(\psi, \nu) = $
$$
 \sup_{p \ge 2} \left[|c|_p(\nu)/\psi(p) \right], \ \ ||c||_m(\nu) =
 \sup_{p \ge 2} \left[|c|_p(\nu) \ \cdot p^{-1/m} \right], \ \ m = \const > 0.
$$

 Note than our  Orlicz $ N \ - $ functions $ N \in EOS $ does not satisfy the so-called
$ \Delta_2 $ condition. \\

\vspace{2mm}
\section{ Boundedness of Hilbert's transform in GLS}
\vspace{2mm}

We consider in this section the case $ T = [-\pi, \pi] $ equipped with the classical
Lebesgue measure, i.e. $ d = 1, $ and the case
of Hilbert's transform in GLS spaces. \par
 Recall that for the integrable function $ f: T \to R $ with the correspondent Fourier
series

$$
f(x) = 0.5 a(0) + \sum_{k=1}^{\infty} \left[a(k) \cos kx + b(k) \sin kx \right],\eqno(3.1)
$$
where as ordinary
$$
a(k) = \pi^{-1} \int_T f(t) \cos(kt) dt, \ b(k) = \pi^{-1} \int_T f(t) \sin(kt) dt,
$$
the Hilbert's transform $ H[f](x) $ may be defined as follows:

$$
H[f](x) = \sum_{n=1}^{\infty} [a(n) \sin(nx) - b(n) \cos(nx)].
$$
 Equivalent definition:

$$
H[f](x) = (2 \pi)^{-1} \ p.v. \ \int_T f(x-t) \cot(t/2) \ dt. \eqno(3.2)
$$
 See in detail, e.g., the classical monograph of A.Zygmund  \cite{Zygmund1}, chapter  11.\par
 Let $ p \in (1,\infty). $ It is known that the operator $ H[f] $ is bounded in all the
 spaces $ L_p = L_p(T). $ The {\it exact value} of the norm

 $$
 K_H(p) \stackrel{def}{=} \sup_{f \in L_p, f \ne 0} |H[f]|_p/|f|_p =
 |H|(L_p \to L_p)
 $$
 was computed by  S.K.Pichorides \cite{Pichorides1}:

$$
K_H(p) = \tan(\pi/(2p)), \ p \in (1,2]; \ K_H(p) = \cot(\pi/(2p)), \ p \in [2,\infty).
\eqno(3.3)
$$
 Let now $ \psi(\cdot) \in \Psi(1,\infty), $ i.e. $ \supp \psi \subset (1,\infty). $
 We define the new $ \psi $ function $ \psi^{(H)}(p), \ p \in (1,\infty) $ as follows:

 $$
 \psi^{(H)}(p) = K_H(p) \cdot \psi(p), \ p \in (1, \infty).\eqno(3.4)
 $$

 {\bf Theorem H.}

 $$
 ||H[f]||G\psi^{(H)} \le 1 \cdot ||f||G\psi, \eqno(3.5)
 $$
where the constant "1" is the  best possible. \par
{\bf Proof  of the upper bound} is very simple. Let $ ||f||G\psi < \infty, $ since
in other case is nothing to prove.  Moreover, we can and will suppose
$ ||f||G\psi = 1, $ or following

$$
\forall p \in (1,\infty) \ \Rightarrow |f|_p \le \psi(p).
$$

 It follows from the Pichorides inequality

 $$
 |H[f]|_p  \le K_H(p) \cdot |f|_p \le K_H(p) \cdot \psi(p) = \psi^{(H)}(p),
 $$
therefore

$$
||H[f]||G\psi^{(H)} \le 1 = ||f||G\psi.
$$

{\bf Proof  of the exactness.} We will use the method of the proposition 1. Let us
denote

$$
V(f,\psi) = \frac{||H[f]||G\psi^{(H)}}{||f||G\psi} =
\frac{\sup_{p > 1} [|H[f]|_p/(K_H(p)\cdot \psi(p))]}{\sup_{p> 1} [|f|_p/\psi(p)]},
$$

$$
\overline{V}= \sup_{\psi \in \Psi} \sup_{f \ne 0, f \in G\psi} V(f,\psi).
$$
 The assertion of theorem H may be formulated as equality $ \overline{V}=1; $ we
proved $ \overline{V}\le 1. $ It remains to prove that $ \overline{V}\ge 1. $ \par
 If we implement the natural choice of the function $ \psi(p) $ for the
 $ f(\cdot): \ \psi(p) = |f|_p, $ we receive the inequality

 $$
 \overline{V}\ge \sup_{p > 1} \sup_{f \in L_p}
 \left[ \frac{|H[f]|_p}{K_H(p)\cdot |f|_p} \right] \ge
\overline{\lim}_{p \to \infty} \sup_{f \in L_p}
\left[ \frac{|H[f]|_p}{K_H(p)\cdot |f|_p} \right].
 $$
Let us consider the {\it family} of a functions

$$
g_{\Delta}(x) = \sum_{n=1}^{\infty} n^{-1} \ \log^{\Delta}(n) \ \sin nx, \
\Delta = \const > 0.
$$
 It is known, see, e.g., \cite{Zygmund1}, chapter 8, that as $ x \to 0 $

 $$
 g_{\Delta}(x) \sim \frac{2}{\pi} |\log |x||^{\Delta},
 $$
therefore as $ p \to \infty $
$$
|g_{\Delta}(\cdot)|_p \sim 2^{1/p} \ \frac{2}{\pi} \
\left[\int_0^{\pi} |\log x|^{\Delta} dx   \right]^{1/p} \sim
$$

$$
2^{1/p} \ \frac{2}{\pi} \ (\Gamma(\Delta p + 1))^{1/p} \sim 2^{1/p} \ \frac{2}{\pi} \
\left[\frac{\Delta p}{e} \right]^{\Delta}.
$$
 Further,

$$
f_{\Delta}(x):= H[g_{\Delta}](x) = \sum_{n=1}^{\infty} n^{-1} \
 \log n^{\Delta} \ \cos(nx),
$$
then as $ x \to 0 $

$$
|f_{\Delta}(x)| \sim \frac{|\log |x||^{\Delta}}{\Delta+1},
$$
and correspondingly as $ p \to \infty, $ i.e. the critical point $ p_0 = \infty:$
$$
|f_{\Delta}(\cdot)|_p \sim 2^{1/p} \ (\Delta+1)^{-1} \
\left[ \frac{(\Delta+1) p}{e}\right]^{\Delta+1};
$$

$$
\frac {|f_{\Delta}(\cdot)|_p}{|g_{\Delta}(\cdot)|_p} \sim \frac{2}{\pi} \cdot
p \ e^{-1} \ \left(1 + \frac{1}{\Delta} \right)^{\Delta}.
$$
 It follows from the Pichorides result that as $ p \to \infty $

$$
K_H(p) \sim \frac{2}{\pi} \ p.
$$
We find substituting into the expression for $ \overline{V} $ for all the values
$ \Delta  > 0: $

$$
 \overline{V} \ge e^{-1} \ \left(1 + \frac{1}{\Delta} \right)^{\Delta}.
$$
The expression in the right hand side tends to one as $ \Delta \to \infty. $ \par
 This completes  the proof of our theorem. \par

\vspace{2mm}

 Analogous result is true for  "continuous"   Hilbert's transform, i.e. in the space
$ X = R^1. $ Recall that  in this case

$$
H[f](x) = \pi^{-1} \ p.v. \ \int_{-\infty}^{\infty} \frac{f(t) \ dt}{x-t}.\eqno(3.6)
$$

Since

$$
K_H(p) \stackrel{def}{=} |H|(L_p \to L_p) \sim \frac{2}{\pi} \ \frac{1}{p-1}, \
p \to 1+0,
$$
the correspondent example may be constructed as follows:

$$
f_0(x) = 1, \ x \in (0,1), \ f_0(x) = 0, \ x  < 0, \ x > 1;
$$
then

$$
H[f_0](x) = \pi^{-1} \log \left|\frac{x}{x-1}  \right|;
$$

$$
x \to \infty \ \Rightarrow \pi H[f_0](x) \sim 1/x;
$$

$$
|\pi H[f_0]|_p \sim  \frac{2}{\pi} \ \frac{1}{p-1} = \frac{2}{\pi} \ \frac{1}{p-1} |f_0|_p, \ p \to 1+0. \eqno(3.7)
$$

 See in detail \cite{Bennet1}, p. 126-128.\par

\vspace{2mm}
\section{ Weight Fourier's  inequalities for GLS spaces}
\vspace{2mm}

Let again $ d = 1, \ X = [-\pi, \pi], $

$$
f(x) = 0.5 a(0) + \sum_{k=1}^{\infty} [a(k) \cos kx + b(k) \sin kx],
$$
$$
\gamma = \const \in (0,1),  U_{\gamma}[f] = |x|^{-\gamma} f(x),
$$
and we define the sequence $ \lambda(n), n=1,2,\ldots $ as follows:
$$
\lambda(1) = a(0), \lambda(2) = b(1), \lambda(3) = a(1), \lambda(4) = b(2),
\lambda(5) = a(2), \ldots.
$$
We intend to obtain in this section the GLS norm estimation for the function
$ U_{\gamma}[f] $ through the GLS norm estimation for the coefficients
$ \{\lambda(n)\}.$\par

\vspace{2mm}

We consider in this section that both the sequences $ a(n) $ and $ b(n) $
are monotonically decreasing; more general case may be investigated by means of the
main result of the article \cite{Yu1}, see also \cite{Beckner1}.\par

 A new notations: $ p_0 = 1/\gamma, $ (critical point);

$$
|\vec{\lambda}|^{(\gamma)}_p =
\left[\sum_{n=1}^{\infty} n^{p(1+\gamma) - 2} |\lambda(n)|^p \right]^{1/p}, \eqno(4.1)
$$

$$
l^{(\gamma)}_p = \{\vec{\lambda}: \ |\vec{\lambda}|^{(\gamma)}_p < \infty \},\eqno(4.2)
$$
and we define for arbitrary function $ \psi \in G\Psi(1,p_0) $

$$
||\vec{\lambda}||G^{(\gamma)}\psi = \sup_{p \in (1,p_0)}
|\vec{\lambda}|^{(\gamma)}_p /\psi(p), \eqno(4.3)
$$

$$
\psi^{(\gamma)}(p) = K^{(\gamma)}(p) \cdot
 \psi(p), \ p \in (1, p_0); \eqno(4.4)
$$
where

$$
K^{(\gamma)}(p) \stackrel{def}{=} \sup_{\lambda \in l^{(\gamma)}_p, \lambda \ne 0}
\left[ \frac{|U_{\gamma}[f]|_p}{|\vec{\lambda}|^{(\gamma)}_p}\right] < \infty, \
 p \in (1, p_0) .\eqno(4.5)
$$

{\bf Theorem $ \gamma.$ }

$$
||f||G \psi^{(\gamma)} \le 1 \cdot ||\vec{\lambda}||G^{(\gamma)}\psi,
$$
where the constant "1" is the best possible. \par
{\bf Proof.} \par
{\bf 1.} It follows  after  some calculations in the article \cite{Yu1}, see also
\cite{Boas1}, \cite{Chen1}, \cite{Chen2} that

$$
K^{(\gamma)}(p) \le \frac{C_{\gamma}(p) \ \gamma^{-2}}{1/\gamma - p},\eqno(4.6)
$$
where $ C_{\gamma}(p) $ is continuous function of the variable $ p $
on the {\it closed} interval $ p \in [1,p_0], $ and such that

$$
\lim_{p \to p_0-0} C(p) = 1.
$$
{\bf 2.} Let us estimate from below the constant $ K^{(\gamma)}(p).$ It is enough to
consider the following example.

$$
g_{\Delta}(x) = \sum_{n=2}^{\infty} n^{-1} \ \log^{\Delta}(n) \ \cos(nx); \
\Delta = \const > 0,\eqno(4.7)
$$
i.e.here $ a(n) = 0, \ n= 0,1,2,\ldots; \ \lambda(n) = b(n) = n^{-1}\log^{\Delta}( n) .$ \par
 It is easy to calculate as $ p \to p_0-0 :$

 $$
\sum_{n=2}^{\infty} \left[n^{p\gamma - 2} \frac{\log^{p \Delta}(n)}{n^p} \right]^{1/p} \sim \frac{\Gamma^{\gamma}(\Delta/\gamma + 1)}{(1-p\gamma)^{\Delta+\gamma}};
 $$

$$
g_{\Delta}(x) \sim (\Delta+1)^{-1} \ |\log |x| \ |^{\Delta+1}, \ x \to 0;
$$

$$
|g_{\Delta}|_p^{(\gamma)} \sim (\Delta+1)^{-1}
\frac{\Gamma^{\gamma}((\Delta+1)/\gamma  )}{(1-p\gamma)^{(1+\Delta)/\gamma + 1}},  \
p  \to p_0 - 0,
$$
and we obtain after dividing as $ p \to p_0-0 $

$$
\overline{K}:= (1-p \gamma) \ \frac{|g_{\Delta}|_p^{(\gamma)}}{|\lambda|_p^{(\gamma)}} \sim \frac{\Gamma^{\gamma}((1+\Delta)+1)}{(\Delta+1) \Gamma^{\gamma}(\Delta/\gamma + 1)}.
$$
 It follows from Stirling's formula that  as $ \Delta \to \infty $

$$
\overline{K} \sim e^{-1} \gamma^{-1} \left[\frac{1+\Delta}{\Delta} \right]^{\Delta} \to
\gamma^{-1}.
$$
So,  one has as $ p \to  p_0-0: $

$$
K^{(\gamma)}(p) \sim \frac{\gamma^{-1}}{1-p \gamma} =
\frac{\gamma^{-2}}{p_0-p}. \eqno(4.7)
$$
 As before, the last assertion proves the proposition of the considered theorem.\par

\vspace{2mm}
\section{ Boundedness of Fourier's  transform}
\vspace{2mm}

{\bf Theorem 1.} {\it Let $ X = [0, 2 \pi]^d  $ and  $ \psi \in \Psi. $ Then the Fourier
operators $ s_M[\cdot] $ are uniformly bounded in the space $ L(N[\psi]) $ into
 another exponential Orlicz's  space} $ L(N_d[\psi]): $

$$
\sup_{M \ge 1} ||s_M[f]||L(N_d[\psi]) \le C_6(d,\psi) \ ||f||L(N[\psi]). \eqno(5.1)
$$

{\bf Theorem 2.} {\it Let now $ X = R^d, \ \psi \in \Psi $ and $ \alpha = \const > 1. $
The Fourier operators $ S_M[\cdot] $ are uniformly bounded in the space $ L(N^{(\alpha)}[\psi]) $
into the space } $ L(N^{(\alpha)}_d[\psi]): $
$$
\sup_{M \ge 1} ||S_M[f]||L(N^{(\alpha)}_d [\psi]) \le C_7(\alpha,d, \psi) \
||f||L(N^{(\alpha)}[\psi]). \eqno(5.2)
$$
 Since the function $ N[\psi] $ does not satisfies the $ \Delta_2 $ condition, the assertions
(5.1) and (5.2) does not mean that in general case when $ f \in L(N_d^{\alpha}[\psi]) $
$$
\lim_{M \to \infty} ||s_M[f] - f||L(N_d[\psi]) = 0, \eqno(5.3)
$$
$$
\lim_{M \to \infty} ||S_M[f] - f||L(N_d^{(\alpha)}[\psi]) = 0; \eqno(5.4)
$$
see examples further. But it is evident that propositions (5.3) and (5.4) are true  if
correspondingly
$$
f \in L^0(N_d[\psi]), \ \ f \in L^0(N_d^{(\alpha)} [\psi]).
$$
Also it is obvious that if $ f \in L(N_d[\psi]), X = [0, 2 \pi]^d $ or, in the case $ X = R^d,\ f \in L(N_d^{(\alpha)}[\psi]), $ then for all EOF $ \Phi(\cdot) $  such that
$ \Phi << N_d[\psi] $ or $ \Phi << N_d^{(\alpha)}([\psi]) $  the following
implications hold:
 $$
 \forall f \in L(N_d[\psi]) \ \Rightarrow \lim_{M \to \infty} ||s_M[f] - f||L(\Phi) = 0, \ X = [0,2\pi]^d; \eqno(5.5)
 $$

$$
\forall f \in L(N_d[\psi]) \ \Rightarrow \lim_{M \to \infty} ||S_M[f]-f||L(\Phi) = 0, \ X = R^d. \eqno(5.6)
$$
{\bf Theorem 3.} {\it Let $ \Phi(\cdot) $ be an EOF and let
$ N(\cdot) =  L^{-1}(u) \in EOF, $ where $ L(y), \  y \ge \exp(2) $ is a positive
slowly varying at $ u \to \infty $
strongly increasing continuous differentiable in the domain $ [\exp(2), \infty) $ function such that the function
$$
W(x) = W_L(x) = \log L^{-1}(\exp x), \ \ x \in [2, \infty)
$$
is again strong increasing to infinity together with the derivative $ dW/dx. $
 In order to the implication (5.5) or,
correspondingly, (5.6) holds, it is necessary and sufficient that
$ \Phi << L(N_d[\psi]), $
or, correspondingly } $ \Phi << L(N_d^{(\alpha)}[\psi]). $ \par
 For instance, the conditions of theorem 3 are satisfied for the functions $ N = N_{m,r}(u). $ \par

\vspace{2mm}

{\bf Theorem 4.} {\it Let} $ \psi \in G\Psi(1,2); $ {\it we denote}

$$
\zeta(p) = \psi(p/(p-1)).
$$
{\it We assert:}

$$
||F[f]||G\zeta \le C(d) \ ||f||G\psi, \ f \in G\psi,
$$
{\it and the last estimation is non-improvable.} \par
{\it As a consequence: let $ f(\cdot) \in
G(1,b, \alpha, 0), \alpha > 0. $ Then $ F[f] \in L(N^{(2)}_{1/\alpha}) $ and }
$$
\sup_{M \ge 1} ||S_M[f]||L \left(N^{(2)}_{1/\alpha} \right) \le C_8(\alpha, N) \
||f||G(1,b, \alpha, 0).
$$

\vspace{2mm}

Analogously may be formulated (and proved) the "discrete" analog of this result.\par

{\bf Theorem 4a.} {\it Let} $ \psi \in G_d\Psi(1,2); $ {\it we denote for the
bilateral complex sequence}

$$
c = \vec{c} = \{\ldots, c(-2), c(-1),c(0),c(1),c(2), \ldots   \}
$$

$$
F[c](x) = F(x) = \sum_{k=-\infty}^{\infty} c(k) \exp(ikx);
$$

$$
\zeta_d(p) = \psi_d(p/(p-1)).
$$
{\it We assert:}

$$
||F[f]||G\zeta_d \le C(d) \ ||c||G\psi_d, \ c \in G_d \psi,
$$
{\it and the last estimation is also non-improvable.} \par

\vspace{2mm}

{\bf Theorem 5 A. } {\it Let $ \{\phi_k(x), k = 1,2, \ldots \} $ be an orthonormal uniform bounded:
$$
M:= \sup_k \vraisup_x |\phi_k(x)| < \infty
$$
sequence of a functions on some non-trivial measurable space $ (X, A, \mu) $ and
 (in the $ L_2(X,\mu) $ sense)}
$$
f(x) = \sum_{k=1}^{\infty} c(k) \ \varphi_k(x). \eqno(5.7)
$$
{\bf A).} {\it \ If $ c \in g(\psi, \nu), $ then }
$$
||f||L \left(N_1[\psi], X, \mu \right) \le C_9 \cdot  (1+  M)
\cdot  ||c||g(\psi,\nu). \eqno(5.8)
$$

 This result may be reformulated as follows. Let $ c = \vec{c} = \{c(k) \}, \
  k = 1,2,\ldots $ be some numerical sequence such that for some
  $ \psi \in G_{d, \nu} \Psi \ c \in  G_{d,\nu}\psi.$ For instance, the function
  $ \psi(p) $ may be natural:  $ \psi(p) = \psi_0(c,p;\nu) $ for the sequence
  $ \{c(k)\} $ relative the $ \nu(\cdot) $
  norm:

  $$
 \psi_0(c,p;\nu):= \left[\sum_{k=1}^{\infty} |c(k)|^p \ k^{p-2} \right]^{1/p},
  $$
 if there exists for some non-trivial interval $ p \in (A,B); 1 \le A < B \le
 \infty. $ \par
  We define also

  $$
  \tilde{\psi} = p \cdot (1+M) \cdot \psi(p).
  $$

{\bf Theorem 5 A'.}

$$
||f||G\tilde{\psi} \le K_3 \ ||c||G_d\psi,
$$
{\it and the last inequality is asymptotically exact.}\par
 The proof is at the same as in the proposition 1.
Note that the point $ p = \infty $ is unique "critical" point in this considerations.\par

\vspace{2mm}

{\bf Theorem 5 B.} {\it Let $ c \in g(\alpha)  $ for some $ \alpha \in (0,1]. $ We assert that }
$$
||f|| L \left(N_{1/\alpha}, X, \mu \right) \le C_{10}(\alpha) \cdot (\max(1, \sup_{k,x} |\phi_k(x)|) \cdot  \ ||c(\cdot)||g(\alpha).
$$

{\bf Theorem 5 C.}
  {\it Let $ \{\phi_k(x), k = 1,2, \ldots \} $ be again some orthonormal uniform bounded:
$$
M:= \sup_k \vraisup_x |\phi_k(x)| < \infty
$$
sequence of a functions on some non-trivial measurable space $ (X, A, \mu) $ and
 (in the $ L_2(X,\mu) $ sense)}
$$
f(x) = \sum_{k=1}^{\infty} c(k) \ \varphi_k(x).
$$
{\it Let } $ c = \vec{c} \in G_d \psi $ {\it for some} $ \psi \in \Psi(1,2). $
{\it Denote }

$$
\tau(q) = (1+M) \cdot \psi(q/(q-1)), \ q \in (2,\infty).
$$
{\it Proposition:}

$$
||c||G_d\tau \le C_5 \ ||f||G \psi,
$$
{\it and the last inequality is asymptotically exact.}\par

\vspace{2mm}

{\bf Theorem 6.} {\it If $ f \in G(\alpha; \psi,\nu), \ $ where $ \alpha \ge 2, $  then }
$$
  \sup_{M \ge 1} ||S_M[f]||L \left(N^{(\alpha)}_d [\psi] \right) \le
  C_{11}(\alpha,\psi, N, \nu) \ ||f||G(\alpha; \psi,\nu). \eqno(5.9)
$$

\vspace{2mm}

\section{ Auxiliary results}

\vspace{2mm}
{\bf Theorem 7.} {\it Let  $ N(u) = N(W,u) = \exp(W(\log u)), \ u > e^2,  \
\psi(p) = \exp(W^*(p)/p), \ p \ge 2, $ and $ X = T^d. $
We propose that the Orlicz's norm $ ||\cdot|| L(N) $ and the  norm
$ ||\cdot||G(\psi) $ are equivalent.  Moreover, in this case  $ f \ne 0, \ f \in G(\psi) $
(or $ f \in L(N(W(\cdot),u) \ ) $ if and only if}  $ \ \exists C_{12},C_{13}, C_{14} \in
(0,\infty) \ \Rightarrow \forall \ u > C_{14} $
$$
T(|f|,u) \le C_{12} \exp \left(-W \left(\log \left(u/C_{13} \right) \right) \right), \eqno(6.1)
$$
{\it where for each measurable function } $ f: X \to R $
$$
T(|f|,u) = \mes \{x: |f(x)| > u \}.
$$

{\bf Proof } of theorem 7. \ A). Assume at first  that $ f \in L(N), \ f \ne 0.$
Without loss of generality we suppose that $ ||f||L(N) = 1/2. $ Then
$$
\int_X N(W, |f(x)|) \ dx \le 1 < \infty.
$$
 The proposition (6.1) follows from Tchebyshev's inequality such that in (6.1)
$ C_{12} = 1, \ C_{13} = C_{14} = 1/||f||L(N), \ f \ne 0. $ \par
B). Inversely, assume that $ f, \ f \ne 0 $ is a measurable function, $ f: X \to R^1 $ such that
$$
T(|f|,u) \le \exp(-W(\log u)), \ u \ge e^2.
$$
 We have by virtue of properties of the function $ W: $
$$
\int_X N(|f(x)|/e^2) \ dx = \int_{ \{x: |f(x) \le e^2\} } + \int_{ \{x: |f(x)| > e^2} =
I_1 + I_2;
$$
$$
I_1 \le \int_X N(1) \ dx = N(1),
$$
$$
I_2 \le \sum_{k=2}^{\infty} \int_{e^k < |f| \le e^{k+1} }
\exp(W(|f(x)|/e^2)) \ dx \le
$$
$$
 \sum_{k=2}^{\infty} \exp((W(k-1)) \ T(|f|,k) \le
  \sum_{k=2}^{\infty} \exp(W(k-1) - W(k)) < \infty.
$$
 Thus, $ f \in L(N(W)) $ and
$$
\int_X N(|f(x)|/e^2) \ dx \le N(1) + \sum_{k=2}^{\infty} \exp(W(k-1) - W(k)) < \infty.
$$
C). Let now $ f \in G(\psi); $ without loss of generality we can assume that $ ||f||G(\psi)=1. $ We deduce for $ p \ge 2:$
$$
\int_X |f(x)|^p \ dx \le \psi^p(p).
$$
We obtain using again the Tchebyshev's  inequality:
$$
T(|f|,u) \le u^{-p} \psi^p(p) = \exp \left[-p \log u + p \log \psi(p) \right],
$$
and  after the minimization over $  p: \ u \ge \exp(2) \ \Rightarrow  $
$$
T(|f|,u) \le \exp \left( - \sup_{p \ge 2} (p \log u - p \log \psi(p)) \right) =
$$
$$
\exp \left( (p \log \psi(p))^*(\log u)  \right) = \exp(-W(\log u)).
$$
D). Suppose now that $ T(|f|,x) \le \exp(-W(\log x)), \ x \ge \exp(2). $
We conclude:
$$
\int_X |f(x)|^p \ dx = p \int_o^{\infty} x^{p-1} T(|f|,x) dx = p \int_0^{\exp(2)} +
$$
$$
p \int_{\exp(2)}^{\infty} \le p \int_0^{\exp(2)} x^{p-1} dx + p \int_{\exp(2)}^{\infty}
x^{p-1} \ T(|f|,x) dx \le
$$

$$
 e^{2p} + \int_{\exp(2)}^{\infty} p x^{p-1} \exp(-W(\log x)) \ dx =
$$
$$
e^{2 p} + p \int_2^{\infty} \exp(py - W(y)) \ dy, \ p \ge 2.
$$
We obtain using Laplace's method and  theorem of Fenchel - Moraux:
$$
\int_X |f(x)|^p dx \le e^{2p} + C^p \exp \left(\sup_{y \ge 2} (py - W(y)) \right) = e^{2p} +
$$

$$
C^p \exp(W^*(p)) = e^{2p} + C^p \exp(p \log \psi(p)) \le C^p \ \psi^p(p).
$$
 Finally, $ ||f||G(\psi) < \infty. $ \par
For example, if $ m > 0, \ r \in R, $ then
$$
 f \in L \left(N_{m,r} \right) \Leftrightarrow  \sup_{p \ge 2} \left[|f|_p \ p^{-1/m} \
\log^{-r} p \right] < \infty \ \Leftrightarrow
$$
$$
T(|f|,u)  \le C_0(m,r) \exp \left( - C(m,r) u^m \left( \log^{-mr} u \right) \right), \ u \ge 2.
$$

{\it Remark 3.} If conversely
$$
T(|f|,x) \ge \exp(-W(\log x)), \ x \ge e^2,
$$
then for sufficiently large values of $ p; \ p \ge p_0 = p_0(W) \ge 2  $
$$
|f|_p \ge C_0(W) \ \psi(p), \ \  C_0(W) \in (0, \infty).\eqno(6.2)
$$

{\it Remark 4.} In this proof we used only the condition $ 0 < \mes(X) < \infty. $ Therefore, our  conclusions in theorem 7 are true in this more general case.\par

{\bf Theorem 8.} {\it Let $ \psi \in \Psi.$ We assert that $ f \in L^0(N[\psi]), $ or,
equally, $ f \in G^0(\psi) $ if and only if }

$$
\lim_{p \to \infty} |f|_p /\psi(p) = 0. \eqno(6.3)
$$
{\bf Proof.} It is sufficient by virtue of theorem 7 to consider only the case of
$ G(\psi) $ spaces.\par
 1. Denote $ G^{00}(\psi) = \{ f: \ \lim_{p \to \infty} |f|/\psi(p) = 0 \}. \ $ Let
$ f \in G^0(\psi), \  f \ne 0. $ Then for arbitrary $ \delta = \const > 0 $ there exists a constant $ B = B(\delta,f(\cdot)) \in (0, \infty) $ such that

$$
||f-f I(|f| \le B) \ ||G(\psi) \le \delta/2.
$$
Since  $ |f| I(|f| \le B)| \le B, $ we deduce
$$
|f I(|f| \le B)|_p/\psi(p) \le B/\psi(p).
$$
We obtain using triangular inequality  for sufficiently large values
$ p: \ p \ge p_0(\delta)=p_0(\delta,B) \ \Rightarrow $
$$
|f|_p /\psi(p) \le \delta/2 + B/\psi(p) \le \delta,
$$
as long as $ \psi(p) \to \infty $ at $ p \to \infty. $ Therefore $ G^0(\psi)
\subset G^{00}(\psi). $\par
(The set $ G^{00}(\psi) $ is a closed subspace of $ G(\psi) $ with respect to the $ G(\psi) $ norm and contains all bounded functions.) \par
2. Inversely, assume that $ f \in G^{00}(\psi). $ We deduce denoting $ f_B = f_B(x) =
f(x) I(|f|> B) $ for some $ B = \const \in (0, \infty):$

$$
\forall Q \ge 2 \ \Rightarrow \lim_{B \to \infty} |f_B|_Q = 0.
$$
 Further,
$$
||f_B||G(\psi) = \sup_{p \ge 2} |f_B|_p/\psi(p) \le \max_{p \le Q} |f_B|_p/\psi(p) +
$$

$$
\sup_{p > Q} |f_B|_p/\psi(p) \stackrel{def}{=} \sigma_1 + \sigma_2;
$$

$$
\sigma_2 = \sup_{p > Q} |f_B|_P/\psi(p) \le \sup_{p \ge Q} (|f|_p/\psi(p)) \le \delta/2
$$
for sufficiently large $ Q $ as long as $ f \in G^{00}(\psi). $  Let us now estimate the
value $ \sigma_1: $
$$
\sigma_1 \le \max_{p \le Q} |f_B|_p/\psi(2) \le \delta/2
$$
for sufficiently large $ B = B(Q). $ Therefore,

$$
\lim_{B \to \infty} || f_B||G(\psi) = 0, \ \  f \in G^0(\psi).
$$
{\bf Theorem 9}. {\it Let $ \psi(\cdot) = \psi_N(\cdot), \theta(\cdot) = \theta_{\Phi}(\cdot) \ $
be a two functions on the classes $ \Psi $ with correspondent $ N \ - $ Orlicz's functions
$ N(\cdot), \Phi(\cdot): $

$$
N(u) = \exp \left \{\left[p \log \psi(p) \right]^*(\log u) \right\},
$$

$$
\Phi(u) = \exp \left \{ \left[ p \log \theta(p) \right]^*(\log u \right \}, \ \ u \ge \exp(2).
$$
We assert that  $ \lim_{p \to \infty} \psi(p)/\theta(p) = 0 $ if and only if $ N(\cdot) >> \Phi(\cdot).$ } \par
{\bf Proof } of theorem 9. A). Assume at first that $ \lim_{p \to \infty} \psi(p)/\theta(p) = 0. $
Denote $ \epsilon(p) = \psi(p)/\theta(p), $ then $ \epsilon(p) \to 0, \ p \to \infty. $ \par
 Let $ \{ f_{\zeta}, \ \zeta \in Z \} $ be arbitrary bounded in the $ G(\psi) $ sense set of a functions:
$$
\sup_{\zeta \in Z} ||f_{\zeta}||G(\psi) = \sup_{\zeta \in Z} \sup_{p \ge 2} |f_{\zeta}|_p/\psi(p)= C < \infty,
$$
then
$$
\sup_{\zeta \in Z} |f_{\zeta}|_p/\theta(p) \le C \epsilon(p) \to 0, \ p \to \infty.
$$
 It follows from previous theorem that $ \forall \zeta \in Z \ f_{\zeta} \in G^0(\theta) $ and that the family $ \{f_{\zeta}, \ \zeta \in Z \} $ has uniform absolute continuous norm. Our assertion follows from lemma 13.3 in the book [6].\par
B). Inverse, let $ \Phi(\cdot) << N(\cdot). $ Let us introduce the measurable function
$ f: X \to R $ such that $ \ \forall x \ge \exp(2) $

$$
\exp \left(-2 \left[p \log \psi(p) \right]^*(\log x) \right) \le T(|f|,x) \le
$$

$$
\exp \left( - \left[p \log \psi(p) \right]^*(\log x) \right).
$$
Then (see theorem 7)
$$
f(\cdot) \in G(\psi), \ \ \ C_{15}(\psi) \ \psi(p) \le |f|_p \le C_{14}(\psi) \ \psi(p), \ \  p \ge 2.
$$
 Since $ f \in G(\psi), \  \Phi << N, $ we deduce that $ f \in G^0(\theta), $
and, following,
$$
\lim_{p \to \infty} |f|_p/\theta(p) = 0.
$$
Therefore, $ \lim_{p \to \infty} \psi(p)/\theta(p) = 0. $ \par
{\bf Theorem 10}. {\it Let now $ X = R^d $ and $ \psi \in \Psi. $
We assert that the norms
$ ||\cdot|| L(N^{(\alpha)}, [\psi]) $ and $||\cdot|| G(\alpha, \psi), \ \alpha \ge 1 $
are equivalent.} \par
{\bf Proof.}
1. Let $ \forall p \ge \alpha \ \Rightarrow |f|_p \le \psi(p), \ f \ne 0. $  From
Tchebychev's inequality follows that
$$
\lim_{v \to \infty} T(|f|,v) = 0.
$$

Let us consider for some sufficiently small value $ \epsilon \in (0, \epsilon_0), \
\epsilon_0 \in (0,1) $ the following integral:
$$
I_{\alpha, N}(f) = \int_X N^{(\alpha)}(\epsilon |f(x)| \ dx = I_1 + I_2,
$$
where
$$
I_1 = \int_{\{x: |f(x)| \le v \} } N^{(\alpha)}(\epsilon |f(x)|) \ dx,
\ I_2 = \int_{\{x: |f(x)|> v \} } N^{(\alpha)}(\epsilon |f(x)|) \ dx.
$$
Since for $ z \ge v $
$$
N^{(\alpha)}(z) \le C_{19}(\alpha, N(\cdot)) \ \cdot N( z),
$$
we have for the set $ X(v) = \{x, |f(x)| > v \} $ and using the result of theorem 7 for the space with finite measure:
$$
I_2 = \int_{X(v)} N^{(\alpha)}(\epsilon |f(x)| ) \ dx \le C_{20}(\alpha,N,\epsilon) \
||f||L(N^{(\alpha)}, X(v)) \le
$$

$$
C_{21} (\alpha, \epsilon, \psi) \sup_{p \ge \alpha}\left[ ||f||L_p(X(v))/\psi(p) \right] \le
C_{21} \sup_{p \ge \alpha} |f|_p/\psi(p) < \infty.
$$
Further, since for $ z \in (0,v) \ \Rightarrow $
$$
N^{(\alpha)}(\epsilon z) \le C_{22}(v,\alpha, \epsilon) \ |z|^{\alpha},
$$
we have:
$$
I_1 \le C_{22}(\cdot) \int_X |f(x)|^{\alpha} \ dx = < \infty.
$$
Thus, $ f \in L(N^{(\alpha)}[\psi]), \ ||f||L(N^{(\alpha)} [\psi]) < \infty.$ \par
2). We prove now the inverse inclusion. Let $ f \in L \left(N^{(\alpha)}[\psi] \right) $ and
$$
||f||L \left(N^{(\alpha)}[\psi] \right) = 1.
$$
Hence for some $ \epsilon > 0 $
$$
\int_X N^{(\alpha)}(\epsilon |f(x)|) dx < \infty.
$$
 It follows from the proof of
 theorem 7 and the consideration of two cases: $ |z| \le v; \ |z| > v $ the
following elementary inequality: at $ p \ge \alpha $ and for all $ z > 0 \ \Rightarrow $
$$
|z|^p \le C_{23} (\alpha, \epsilon,N) \ N^{(\alpha)} (\epsilon |z|) \cdot \psi^p(p).
$$
 We obtain for all values $ p, \ p \ge \alpha: $
$$
\int_{R^d} |f(x)|^p \ dx \le C_{24}^p(\alpha, \epsilon,\psi) \ \psi^p(p), \ \ \ ||f||G(\alpha;\psi) <
\infty.
$$

\vspace{2mm}
\section{ Proofs of main results. }
\vspace{2mm}
 At first we consider the case Orlicz spaces, i.e. if the function $ f $ belongs to some
exponential Orlicz space.\par
{\bf Proof of theorems 1,2.} Let $ X = [0, \ 2\pi]^d $ and $ f \in L(N[\psi]) $ for some
$ \psi \in \Psi.$ Without loss of generality we can assume that $ ||f||L(N[\psi]) = 1. $
From theorem 7 follows that
$$
\forall p \ge 2 \ \Rightarrow \ |f|_p \le C_{25}(\psi) \ \psi(p).
$$
 From the classical theorem of M. Riesz follows the inequality:
$$
\sup_{M \ge 1} |s_M[f]|_p \le K_1^d \ p^d \ \psi(p), \ \ K_1 = 2 \pi.
$$
It follows again from theorem 7  that
$$
\sup_{M \ge 1} ||s_M[f] - f||L(N_d[\psi]) \le K_1^d + 1 < \infty.
$$
 For example, if $ N(u) = N_m(u) = \exp(|u|^m)-1 $ for some $ m = \const \ge 1, $ then
$$
\sup_{M \ge 1} ||s_M[f] -f ||L(N_{m/(dm+1)}) \le C_{26}(d,m) \ ||f||L(N_m).  \eqno(7.1)
$$
 The "continual" analog of M.Riesz's inequality, namely, the case $ X = R, \ L(N) =
L_p(R), \ p \ge 2: $
$$
\sup_{M \ge 2} |S_M[f]|_p \le K_2^d \ p^d \ |f|_p, \ \ K_2 = 1
$$
 is proved, for example, in \cite{Wolff1}, p.187 - 188. \par
 This fact permit  us to prove also theorem 2. \par
{\bf Lemma 1.} {\it We assert that the "constant" $ m/(dm+1) $ in the estimation (3.5)
is exact. In detail, for all $ m \ge 1 $ there exists $ g = g_m(\cdot) \in L(N_m) $ such that
$ \forall  \Delta \in (0,1/2) $ }
$$
\sup_{M \ge 1} ||s_M[g]||L \left(N_{(m-\Delta)/(dm+1)} \right) = \infty.
$$
{\bf Proof} of lemma 1. It is enough to prove that
$$
\exists g \in L(N_m), \ \ ||H[g]||L(N_{(m-\Delta)/(dm+1)}) = \infty,
$$
where $ H[g] $ denotes the Hilbert transform on the $ [0, \ 2 \pi]^d, $
see \cite{Edwards1}, p. 193 - 197. Also it is enough to consider the case
$ d = 1. $\par
 Let us introduce the function
$$
g(x) = g_m(x) = | \log(x/(2 \pi))|^{1/m}.
$$
Since for $ u > 0 $
$$
\mes \{x: g_m(x) > u \} = \exp \left(- u^m \right),
$$
we conclude $ g_m(\cdot) \in L(N_m) \setminus L^0(N_m) $ (theorem 7).
 Further, it is very simple to verify using the formula for  Hilbert transform  that
$$
C_{28}(m) \left(| \log (x/ (2 \pi))|^{(m+1)/m} + 1 \right) \le |H[g_m](x)| \le
$$
$$
C_{29}(m) \left( |\log (x/(2 \pi))|^{(m+1)/m} + 1 \right).
$$
 Hence  $ \ \forall u \ge 2 $
$$
\exp \left( -C_{29}(m) u^{m/(m+1)} \right) \le \mes \{ x, |H[g_m](x) > u \} \le
$$

$$
\exp \left(-C_{30}(m) \ u^{m/(m+1)} \right).
$$
 It follows again from theorem 7 that
$$
H[g_m] \in L(N_{m/(m+1)}) \setminus L(N_{m/(m+1)}).
$$
Thus $ \forall \Delta \in (0,) \ \Rightarrow \ H[g_m]
\notin L(N_{(m-\Delta)/(m+1)}). $ \par
{\bf Proof} of theorem 3. Let us consider the following function:
$$
z(x) = z_L(x) = \sum_{n=8}^{\infty} n^{-1} \ L(n) \ \sin(nx). \eqno(7.2)
$$
 It is known from the properties of slowly varying functions ([14], p. 98 - 101) that
the series (7.2) converge a.e. and at $ x \in (0, \ 2 \pi] $
$$
C_0 L(1/x) \le z(x) \le C L(1/x).
$$
Therefore, at $ u \in [\exp (2), \infty) $
$$
  L^{-1}(C u)  \le  T(|z|,u) \le   L^{-1}(C_0 u).
$$
 It follows from theorem 7 and (5.2)  that
$$
z(\cdot) \in L(N) \setminus L^0(N), \ N(u) = L^{-1}(u), \ u \ge \exp(2).
$$
 From theorem 8 follows that
$$
0 < C_0 \le |z|_p/\psi(p) \le C < \infty, \ p \ge 2, \ \psi(p) = \exp(W^*(p)/p). \eqno(7.3)
$$
 Note as a consequence that the series (7.2) {\it does not converge in the $ L(N) $ norm,} as long as the system of functions $ \{\sin(nx) \} $ is bounded and hence in the case
when the series (7.2) {\it converge}  in the $ L(N) $ norm $ \Rightarrow z(\cdot) \in L^0(N). $   \par
 Let us suppose now that for some EOF $ \Phi(\cdot) $ with correspondence function $ \theta(p) $ (7.2) convergence in the $ L(\Phi) $ norm. Assume converse to the
 assertion of theorem 3, or equally that
$$
\overline{\lim}_{p \to \infty} \theta(p)/\psi(p) > 0. \eqno(7.4)
$$
 Since the system of functions $ \{ \sin(nx) \} $ is bounded, $ z(\cdot) \in L^0(\Phi). $
By virtue of theorem 8 we conclude that
$$
\lim_{p \to \infty} |z|_p/\theta(p) = 0.
$$
 Thus, we obtain from (7.3)
$$
\lim_{p \to \infty} \psi(p)/\theta(p) = 0,
$$
in contradiction with (7.4). The cases $ X = [0, 2 \pi]^d, X = R^d $ are considered as well as the case $ X = [0, 2 \pi]. $ \par
 Now we consider the case when $ f \in G(a,b,\alpha, \beta).$ \par
{\bf Proof } of theorem 4. Let $ f \in G\psi, \ \psi \in G\Psi(1,2), \
||f||G\psi = 1; $ then

$$
|f|_q \le \psi(q), \ q \in (1,2).
$$

 We denote $ p = q/(q-1), $ then $ p \in [2, \infty). $ We will use the classical result of Hardy - Littlewood, Hausdorff - Young \cite{Pinsky1}, p.193; \cite{Wong1}, p. 93:
$$
|F[f]|_p \le C(d) \ |f|_q \le C(d) \psi(q) = C(d) \psi(p/(p-1)) = C(d) \zeta(p),
$$
or equally
$$
||F[f]||G\zeta \le   C(d) = C(d) \ ||f||G\psi.
$$
 In order to prove the exactness of theorem 4, we consider the following example. Let
 $ d = 1, \ X = R^1, $

$$
f_0(x) = f(x) = |x|^{-1} \ I(|x| \ge 1).
$$
 We deduce for the values $ q \in (1,2) $  and following $ p \in (2,\infty); \
 q \to 1+0 \ \Rightarrow p \to \infty: $

 $$
 |f|_q^q  = 2 \int_1^{\infty} x^{-q} \ dx =  2 \ (q-1)^{-1},
 $$

$$
q \to 1+0 \Rightarrow |f|_q \sim 2 \ (q-1)^{-1};
$$

$$
F[f](t):= F(t) = 2 \int_1^{\infty} x^{-1} \ \cos(tx) \ dx.
$$
 Note that as $ t \to 0+ $

$$
F(t) = 2 \int_t^{\infty} y^{-1} \ \cos(y) \ dy \sim 2 \int_t^1 y^{-1} \ \cos(y) \ dy
\sim
$$

$$
2 \int_t^1 dy/y = 2 |\log t|.
$$
 Therefore, as $ p \to \infty $

$$
|F|_p^p \sim 2^p \int_0^1 |\log t|^p \ dt  = 2^p \ \Gamma(p+1),
$$

$$
|F|_p \sim = \Gamma^{1/p}(p+1) \sim C_2 \ p, C_2 = e^{-1}.
$$
 We conclude after dividing:

$$
\overline{\lim}_{p \to \infty} \frac{F[f_0]_p}{|f_0|_q} \ge C_2/2 > 0,
$$

The second proposition of this theorem follows  from theorem 10.\par

{\bf Proof } of theorem 4a may be ground alike the proof of
of theorem 4 with the analogous counterexample; the Hardy-Young inequality for the
Fourier series has a view

$$
|F[c]|_p \le C_1(d) \  |c|_{q,d}, \ q \in (1,2), \ p = q/(q-1).
$$

\vspace{2mm}

 {\bf Proof of theorem 5A. Upper estimation.} We will use the classical result
of  W.Paley and F.Riesz ( \cite{Wong1}, p.120). \par
Let $ \{ \phi_k(x), \ k = 1,2,\ldots \} $ be some orthonormal bounded sequence of functions. Then $ p \ge 2 \ \Rightarrow $
$$
|f|_p \le K_3 \cdot p \ \cdot \left( 1 + \sup_{k,x} |\phi_k(x)| \right) \cdot
 \left( \sum_k |c(k)|^p \left( |k|^{p-2} + 1 \right) \right)^{1/p}, \eqno(7.5)
$$
where $ K_3 $ is an absolute constant, $ f(x) = \sum_k c(k) \ \phi_k(x). $ \par
Let $ c(\cdot) \in g(\psi, \nu) $ and $ ||c||g(\psi,\nu) = 1. $
By definition of the $ ||\cdot||g(\psi,\nu) $ norm
$$
  \sum_{k = 1}^{\infty} |c(k)|^p \ \left(k^{p-2} + 1 \right)  \le \ ||c||^p
  g(\psi, \nu) \cdot \psi^p(p).
$$
Therefore
$$
|f|_p \le  K_3 \ C_{31} \ p \cdot \psi(p),
$$
and by virtue of theorem 7 $ \ f(\cdot) \in L(N_1[\psi]). $ \par
\vspace{2mm}
{\bf Proof of theorem A.  Exactness.}\par
 Let us consider the following example:

$$
g(x) = \sum_{n=2}^{\infty} n^{-1} \ \log^m n \ \cos(nx), \ x \in [-\pi,\pi].
$$
$ m = \const > 0. $  Here

$$
c(n) = n^{-1} \ \log^m n,
$$
and following as $ p \to \infty $

$$
|c|^p_{p,\nu} = \sum_{n=2}^{\infty} h^{-2} \ \log^{pm} n \ \sim
\int_1^{\infty} x^{-2} \ \log^{pm} x \ dx = \Gamma(pm+1);
$$

$$
|c|_{p,\nu} \sim [\Gamma(pm+1)]^{1/p} \sim m^m e^{-m} \ p^m.
$$
 We know that

 $$
|g(x)| \sim C(m) \ |\log |x| \ |^{m+1}, x \to 0;  \ |g|_p \sim C_1(m) \ p^{m+1}.
 $$
  Substituting into the expression for the value $ \overline{V}, $ we get to the
  assertion of our theorem. \par

\vspace{2mm}

{\bf Proof of Theorem 5 B.} Here we use the "discrete" inequality of Hausdorff-Young,  Hardy - Littlewood (see \cite{Edwards1}, p.101; \  \cite{Gord1}, \cite{Krein1},
chapter 5, \cite{Titchmarsh1}, chapter 4, sections 1,2:
$$
|f|_p \le K_4 \ |c|_q, \ p \ge 2, \ q = p/(p-1), \ \ K_4 = 2 \pi.
$$
If $ ||c||g(\alpha) =1, $ then
$$
|c|_q \le (q -1)^{\alpha},    |f|_p \le K_4 \ p^{\alpha}, \ p \ge 2.
$$
Again from theorem 7 follows that $ f \in L \left(N_{1/\alpha} \right). $ \par
\vspace{2mm}
{\bf Proof of Theorem 5 C} is at the same as the proof of theorem 5. We use
the following classical inequality:

$$
|f|_{p/(p-1)} \le (1+M) \ |c|_p,
$$
see \cite{Kaczmarz1}, chapter 6, section 3.\par
\vspace{2mm}
{\bf Proof } of theorem 6. The analog of inequality (7.5) in the case
$$
F[f](t) = \int_R \exp(i t x) f(x) dx, \ d = 1,
$$
namely:
$$
|F[f]|_p \le K_5 \ p \ |f|_p(\nu),
$$
when $ f(\cdot) \in L_p(\nu) \subset G(\nu, \alpha, \psi), $ see, for example, in
\cite{Novikov1}, p. 108.  Hence, for all $ \ p \ge \alpha $
$$
||F[f]||L \left(N_d^{(\alpha)}[\psi] \right) \le K_5 \ ||f||G(\alpha; \psi, \nu).
$$
(The generalization on the case $ d \ge 2 $ is evident). \par
 Note that the moment estimations for the wavelet transforms and Haar series are described for example in the books
 \cite{Edwards1}, p.21,  \cite{Pinsky1},  p.297  etc.\par
    It is easy to generalize our results on the cases Haar's or wavelet series and transforms.\par
 In detail, it is true in this cases the moment estimation for the partial sums (wavelet's or Haar's)
$$
|P_M[f]|_p \le K_6 \ |f|_p, \ X = [0,1], \ p \ge 1,
$$
where $ K_6 = 13 $ for Haar series on the interval $ X = [0,1] $  and $ K_6 = 1 $ for the classical
wavelet series, in both the cases $ X = [0,1] $ and $ \ X = R.$
 Hence $ \forall \psi \in \Psi, \ f(\cdot) \in L(N[\psi]) $
$$
\sup_{M \ge 1} ||P_M[f] - f||L(N[\psi]) \le (K_6+1) \ ||f||L(N[\psi]). \eqno(7.6)
$$
 But (4.6) does not mean in general case the convergence
$$
\lim_{M \to \infty} ||P_M[f] - f||L(N[\psi]) = 0, \eqno(7.7)
$$
as long as if (7.7) is true, then $ f(\cdot) \in L^0(N[\psi]) $ and conversely if
$ f \in L^0(N[\psi]), $ then (7.7) holds.\par
 For the different generalizations of wavelet series the estimation (7.6) with constants
$ K_6 $  not depending on $ p, \ p \ge 2 $ see, for example, in the books
\cite{Edwards1}, \cite{Pinsky1}, \cite{Wong1}  etc.\par

\vspace{2mm}
\section{ Concluding remarks. Maximal operators.}
\vspace{2mm}

 We consider in this section the so-called {\it maximal Fourier operators} and
  investigate their boundedness in some Grand Lebesgue spaces. \par
  Let us define the following maximal operators in the one-dimensional case $ d = 1:$
  
  $$
  s^*[f](x) = \sup_{M \ge 1} |s_M[f](x)|, \eqno(8.1),
  $$
  
$$
F^*[f](x) = \sup_{a > 0} \left|\int_{-a}^a f(t) \exp(i t x) \ dt  \right|, \eqno(8.2),
$$  
  
$$
R^*[f](x) = \sup_{a > 0} 
\left| \int_{-\infty}^{\infty} f(t) \frac{\sin(a(x-t))}{x-t} \ dt \right|. \eqno(8.3)
$$  
  
\vspace{2mm}

{\bf Theorem M1.} Let $ f \in G \psi, \supp \psi = (1,\infty). $   We define 

$$
\psi_{\lambda,\mu}(p) = \frac{p^{\lambda}}{(p-1)^{\mu}} \ \psi(p). 
$$
 We assert:
 
 $$
 ||s^*[f]||G\psi_{4,3} \le K_{4,3} \ ||f||G\psi, \eqno(8.4),
 $$  
  
$$
||R^*[f]||G\psi_{4,2} \le K_{4,2} \ ||f||G\psi. \eqno(8.5)
$$  
  
\vspace{2mm}

{\bf Theorem M2.} Let $ f \in G \psi, \supp \psi = (1,2). $   We define as before 

$$
\zeta(q) =  q^2 \ \psi(q/(q-1)), \ q \in (2,\infty). 
$$
Assertion:

$$
||F^*[f]||G\zeta \le K_5 \ ||f||G\psi. \eqno(8.6)
$$
{\bf Proof} is at the same as before. It used the following maximal $ L_p $ 
Fourier estimations:

$$
|s^*[f]|_p \le K_{4,3} \ |f|_p \ \frac{p^4}{(p-1)^3}, \ p \in (1,\infty);
$$

$$
|R^*[f]|_p \le C \ K_{4,2} \ |f|_p \ \frac{p^4}{(p-1)^2}, \ p \in (1,\infty);
$$

$$
|F^*[f]|_q \le  K_5 \ |f|_p \ (p-1)^{-2}, \ p \in (1,2), \ q = p/(p-1);
$$
see, e.g., the classical monograph of Reyna \cite{Reyna1}, p. 144-152; or 
\cite{Beckner1}. \par

\end{document}